\documentclass[notitlepage,11pt,reqno]{amsart}
\usepackage{amssymb}
\usepackage{amsmath}
\usepackage{xspace}
\usepackage[mathscr]{eucal}

\newtheorem{theorem}{Theorem}
\newtheorem{lemma}{Lemma}
\newtheorem{cor}{Corollary}
\newtheorem{prop}{Proposition}

\newcommand{\toward}[1]{\ensuremath{\overset{\rightarrow}{#1}}\xspace}
\newcommand{\from}[1]{\ensuremath{\overset{\leftarrow}{#1}}\xspace}
\newcommand{\N}{\ensuremath{\mathbb{N}}\xspace}
\newcommand{\Z}{\ensuremath{\mathbb{Z}}\xspace}
\newcommand{\nn}{\ensuremath{\mathbb{N}_0}\xspace}
\newcommand{\ssumm}{\ensuremath{\sum_{j=0}^m s_j}}
\newcommand{\aaa}{\ensuremath{\alpha}}
\newcommand{\grint}[1]{\ensuremath{\lceil #1 \rceil}}
\newcommand{\leint}[1]{\ensuremath{\lfloor #1 \rfloor}}
\newcommand{\fract}[3]{\ensuremath{\frac{{#1}^{#3} - {#2}^{#3}}{{#1} -
{#2}}}}

\begin{document}
\title{Divisibility Tests and Recurring Decimals in Euclidean domains}
\author{Apoorva Khare}
\address{Department of Mathematics (University of Chicago), 5734 S.
University Avenue, Chicago, IL - 60637}
\email{\small{apoorva@math.uchicago.edu}}
\date{\today}

\begin{abstract}
\small{In this article, we try to explain and unify standard divisibility
tests found in various books. We then look at recurring decimals, and
list a few of their properties. We show how to compute the number of
digits in the recurring part of any fraction. Most of these results are
accompanied by a proof (along with the assumptions needed), that works in
a Euclidean domain.

We then ask some questions related to the results, and mention some
similar questions that have been answered. In the final section (written
jointly with P. Moree), some quantitative statements regarding the
asymptotic behaviour of various sets of primes satisfying related
properties, are considered.}\\
\end{abstract}
\maketitle

\tableofcontents

\section*{\bf Part 1 : Divisibility Tests}

\section{The two divisibility tests: going forward and
backward}\label{secKh}


We are all familiar with divisibility tests for a few integers such as 3,
9, and 11. To test whether a number is divisible by 3 or 9, we look at
the sum of its digits, which is equivalent to taking the {\it weighted}
sum of the digits, the weights all being 1. For 11, the corresponding
test is to examine the alternating sum and difference of digits, which
means that the weights are $1, -1, 1, -1, \dots$ .

The problem of finding a sequence of weights for various divisors has
been dealt with in \cite{Kh}; in this section, we briefly mention the
various tests. First, some notation. Any $s \in \N$ with digits $s_j : 0
\leq j \leq m$ is a polynomial of 10, i.e.
\[ s = \underline{s_m s_{m-1} \dots s_1 s_0} = \ssumm 10^j =:
\from{s}(k) \]

\noindent and similarly, one defines $\from{s}(k)$, given $s,k \in \N$.
The coefficients are the {\it digits}, i.e. $0 \leq s_j \leq 9$ for all
$j$. To test whether $s$ is divisible by $d \in \N$, let $k \equiv 10
\mod d,\ |k|$ minimal. Then
\[ s \equiv \from{s}(10) \equiv \from{s}(k) \mod d \]

\noindent and thus $s$ is divisible by $d$ if and only if $\from{s}(k)$
is, i.e. the sequence of weights is simply $k^0, k^1, k^2, \dots$, where
\begin{equation}
k = \from{k}_d \equiv 10 \mod d.
\end{equation}

\noindent (In any base $B > 1$, this test would correspondingly use
$\from{k}_d \equiv B \mod d$.)\\

This proves the divisibility tests for 3, 9, and 11 (the last one holds
because $10 \equiv -1 \mod 11$). It also helps test for powers of 2 and 5
(e.g. a number is divisible by 25 iff its last two digits are 00, 25, 50,
or 75). However, if $n > 20$, then $k \equiv 10 \mod n$, and the one with
smallest $|k|$ is 10 itself. So this does not simplify the testing of
divisibility.\\

We now mention the reverse test, also given in \cite{Kh}. Suppose $d \in
\N$ is now coprime to both 2 and 5. Then 10 has a multiplicative inverse
modulo $d$, and we shall take $k$ to be the representative of least
absolute value. Thus, $10k \equiv 1 \mod d$ and $|k|$ is minimal. Then
given $s = \underline{s_m s_{m-1} \dots s_1 s_0} = \ssumm 10^j$, we have
\begin{equation}\label{E2}
s \equiv \ssumm 10^j \equiv \ssumm k^{-j} \equiv 10^m \bigg( \ssumm
k^{m-j} \bigg) \mod d.
\end{equation}

\noindent Since $gcd(d,10) = 1,\ s$ is divisible by $d$ if and only if
\[ \toward{s}(k) = \ssumm k^{m-j} \]

\noindent is too. Thus, the sequence of weights (proceeding right to left
this time) is simply\\ $k^0, k^1, k^2, \dots$, where
\begin{equation}\label{E4}
k = \toward{k}_d \equiv 10^{-1} \mod d.
\end{equation}

\noindent Observe that we can also use equation \eqref{E2} to compute $s
\mod d$ even when $s$ is not divisible by $d$. We simply reduce $10^m
\mod d$ and multiply the result by $\toward{s}(\toward{k}_d) \mod d$.

For example, $\toward{k}_{29} = 3$. Thus, to test 841 for divisibility by
29, we compute $\overrightarrow{841}(3) = 8 + (4 \times 3) + (1 \times
3^2) = 8+12+9 = 29$, which is divisible by 29. Hence $29|841$. Note that
this test is better than the standard one, because $\from{k}_{29} = 10$
or $-19$. (As a later result on $k$-values for arithmetic progressions
shows, we might, of course, get large values of $k$ for large values of
$d$.) Note that the tests for 3, 9, and 11 work this way just as well.\\

\section{The proofs for a general Euclidean domain}

\subsection{Preliminaries}\hfill\break

\noindent {\bf Definitions.}
\begin{enumerate}
\item Let $(E,\nu)$ be a (commutative) Euclidean domain with unity. Thus,
$E$ is a commutative integral domain, with unity and a {\it valuation}
$\nu : E \setminus \{0\} \to \nn := \N \cup \{0\}$, such that:
\begin{enumerate}
\item $\forall s,0 \neq d$ in $E$, we can carry out the {\it Euclidean
algorithm}. In other words, there exist $q,r \in E$ so that $s = dq+r$,
with $r=0$ or $\nu(r) < \nu(d)$.
\item If $a,b \in E$ are nonzero, then $\nu(ab) \geq \nu(a)$.\\
\end{enumerate}

\noindent (By general theory, $E$ is also a PID, and hence a UFD.)\\

\item Define $E_0 := E^\times \cup \{ 0 \}$, where $E^\times$ is the
group of units in $E$. (Thus, $E_0 = E$ iff $E$ is a field.)\\

\item Given $(s,d)$ as above, define $q_0 = s$. Given $q_n$, inductively
define $q_{n+1}$ and $r_n$ to be any choice of quotient and remainder,
when one divides by $d$. In other words, $q_n = d q_{n+1} + r_n$ for all
$n$.

We say that the Euclidean algorithm {\it terminates} for $(s,d)$, if
there is a sequence of $q_n$'s, that is eventually identically zero
(equivalently by Lemma \ref{L0} below, some $q_n = 0$).\\

\item We also say that the {\it strong Euclidean algorithm holds} if
given $s \in E$ and a nonzero nonunit $d \in E \setminus E_0$, for any
choice of quotient $q$ and remainder $r$ (in the Euclidean algorithm), we
have $\nu(q) < \nu(s)$ (resp. $\nu(r) < \nu(d)$), if $q$ (resp. $r$) is
nonzero.\\
\end{enumerate}

\noindent For example, the strong Euclidean algorithm holds for $E = \Z$
and $F[X],\ F$ a field. The latter is easy to see, since the degrees of
$q,r$ must be less than that of $s,d$ respectively; for the former,
suppose $s=dq+r$. If $|dq| \leq |s|$ then we are done; else $|dq|>|s|$.
But then $r \neq 0$, so $|q|= |s-r|/|d| \leq (|s| + |r|)/|d| < (|s|/|d|)
+ 1$. But this is at least $(|s|/2) + 1$ since $d \notin E_0$, so $|q| <
(|s|/2) + 1$. If $|s| > 2$ then we are done; otherwise $d > 2 = s$, and
we have $2 = d \cdot 1 + (2-d) = d \cdot 0 + 2$. In both cases, $\nu(s) =
2 > \nu(q)$.\\

We now show two easy lemmas; the first shows what it means to satisfy one
of the defining axioms of a Euclidean domain, and the second explores
properties of units, related to the valuation.

\begin{lemma}\label{L0}
Given a commutative integral domain $E$ with unit, that satisfies the
Euclidean algorithm (as in (1)(a) above) for some $\nu : E \to \N_0$, the
following are equivalent:
\begin{enumerate}
\item The Euclidean algorithm gives unique quotient and remainder, when
applied to $(de,d)$ for any $e \in E$ and any $d \neq 0$.
\item The Euclidean algorithm gives unique quotient and remainder, when
applied to $(0,d)$ for any $d \neq 0$.
\item $\nu(ab) \geq \nu(a) \forall a,b \neq 0$.
\end{enumerate}
\end{lemma}

\begin{proof}
Clearly, (1) implies (2). Given (2), suppose $\nu(ab) < \nu(a)$ for some
$a,b$. Then if we divide $0$ by $a$, we have two solutions to the
Euclidean algorithm, namely
\[ 0 = 0 \cdot a + 0 = (-b) \cdot a + (ab) \]

\noindent which is a contradiction. Finally, given (3), suppose we have a
solution to $de = d \cdot q + r$. Thus $r = d(e-q)$, so that if $e \neq
q$, then $\nu(r) \geq \nu(d) > \nu(r)$, a contradiction. Hence $e=q$ and
$r=0$.
\end{proof}\hfill

\begin{lemma}\label{L1}
The following are equivalent for $0 \neq B \in E$:
\begin{enumerate}
\item $B \in E^\times$ (i.e. $B$ is a unit).
\item $\nu(Be) = \nu(e)$ for all nonzero $e \in E$.
\item $\nu(B) = \nu(1)$.
\item $\nu(B) \leq \nu(1)$.
\item The set $\{ 0 \neq a \in E : \nu(a) < \nu(B) \}$ is empty.
\end{enumerate}
\end{lemma}

\begin{proof}
Firstly, if $B$ is a unit, then $\nu(e) \leq \nu(Be) \leq \nu(B B^{-1} e)
= \nu(e)$, so they are all equal and $(1) \Rightarrow (2)$. Setting
$e=1,\ (2) \Rightarrow (3) \Rightarrow (4)$. Next, (4) implies (5)
because $\nu(a) = \nu(a \cdot 1) \geq \nu(1) \geq \nu(B)$ for all nonzero
$a$, so that the desired set of elements is empty. Finally, assume (5).
Then $1 = qB+r$ by the Euclidean algorithm, and $r=0$ or $\nu(r) <
\nu(B)$. Given our assumption (5), we conclude that $r=0$ and $1=qB$, so
that $B$ is a unit.
\end{proof}\hfill

We now adapt the results of the previous section to $(E,\nu)$. Henceforth
we assume that $E$ is not a field. Fix a base $B \notin E_0$. From Lemma
\ref{L1} above, this means that we can talk of ``digits" for numbers. Fix
also a divisor $d \in E$ coprime to $B$ (i.e. g.c.d.$(d,B) = 1$).

\begin{lemma}\label{Lprep2}
\hfill
\begin{enumerate}
\item If the strong Euclidean algorithm holds, then it terminates for all
$(s,B)$. In other words, every $s \in E$ can be written as $s =
\underline{s_m s_{m-1} \dots s_1 s_0} = \ssumm B^j$, where $s_i = 0$ or
$\nu(s_i) < \nu(B)$ for all $i$.
\item If $(B,d) = 1$, then $\exists !\ B^{-1} \mod d$.\\
\end{enumerate}
\end{lemma}

\begin{proof}
The first part follows from carrying out the Euclidean algorithm
repeatedly, i.e. repeatedly dividing by $B$. Note that at each stage, the
quotient $q_{n+1}$ has a strictly smaller $\nu$-value than the previous
quotient $q_n$, whence it must eventually be less than $\nu(B)$. We then
define $q_{n+2} = 0$ and $r_{n+1} = q_{n+1}$.

The second part holds in any PID: since $B$ and $d$ are coprime we can
find $\alpha_0, \beta_0 \in E$ such that $\alpha_0 B + \beta_0 d = 1$. Thus
$\alpha_0 B \equiv 1 \mod d$, and $\alpha_0$ is an inverse modulo $d$.

If there are two such, namely $B_1^{-1}$ and $B_2^{-1}$, then $(B_1^{-1}
- B_2^{-1})B \equiv 1 - 1 \equiv 0 \mod d$, so it equals $ld$ for some $l
\in E$. Multiplying by $\alpha_0$, and noting that $\alpha_0 B = 1 -
\beta_0 d$, we get
$$(B_1^{-1} - B_2^{-1}) = [(B_1^{-1} - B_2^{-1}) \beta_0 + l \alpha_0]d$$

\noindent which means that $d | LHS$, and $B^{-1}$ is unique modulo $d$,
as claimed.
\end{proof}\hfill

\subsection{The proof of the reverse test}\hfill\break

We again test whether or not $s$ is divisible by $d$ in $E$. Clearly, we
have two different methods, as above. The first (or standard) method is
as above, and we focus on the second one now. For notation we fix $B >
1$, our base, and $d$, our divisor, that is coprime to $B$. Given a
dividend $s \in E$, let $s_0$ be the unit's digit and $\beta$ be the
``rest", i.e. $s = \beta B + s_0$, with $s_0 = 0$ or $\nu(s_0) <
\nu(B)$.\\

\begin{theorem}\label{Tprep1}
There is a unique $k \in E / d E$, that satisfies the following condition
for all $s \in E$ :
$$s = \underline{\beta s_0} \equiv 0 \mod d \text{ iff } \beta + k s_0
\equiv 0 \mod d$$

\noindent Moreover,
\begin{equation}
k = k_d := \toward{k_d} \equiv B^{-1} \mod d
\end{equation}
\end{theorem}

\begin{proof}
Existence is easy: define $k_d \equiv B^{-1} \mod d$. Now we see that
$\beta + k_d s_0 \equiv k_d(B \beta + s_0) \equiv k_d s \mod d$. Since
$B$ and $d$ were coprime, hence so are $k_d$ and $d$, so $s = B \beta +
s_0 \equiv 0 \mod d$ iff $\beta + k_d s_0 \equiv 0 \mod d$, as desired.

For uniqueness, suppose any such $k$ exists. We then keep $k_d$ as above.
Now take any multiple $s$ of $d$, so that $s_0 = 1$, i.e. $s \equiv 1
\mod B$. Since $d$ and $B$ are coprime, such a multiple clearly exists
(e.g. from the proof of Lemma \ref{Lprep2} above, $\beta_0 d \equiv 1
\mod B$). Then by our condition, $\beta + k s_0 = \beta + k \equiv 0 \mod
d$. But since $s \equiv 0 \mod d$, hence $k_d s \equiv 0 \mod d$, whence
we get $k_d s = k_d (B \beta + 1) \equiv \beta + k_d \equiv 0 \mod d$.

Thus $\beta + k_d \equiv \beta + k \equiv 0 \mod d$, whence $k \equiv
k_d \mod d$, as required.
\end{proof}\hfill

\noindent{\bf Notation :} 
Henceforth, we only consider dividends of the form $s = \underline{s_m
\dots s_0}$, where the $s_i$'s are the digits of $s$, when written in
base $B$. Let $0 < n_1 < \dots < n_l < m$ for some $l$. Then define
$\aaa_{n_1,\dots,n_l,s}(k_d)$ to be the sum
$$\underline{s_m \dots s_{n_l}} + k_d^{n_l - n_{l-1}} \underline{s_{n_l
-1} \dots s_{n_{l-1}}} + k_d^{n_l - n_{l-2}} \underline{s_{n_{l-1} -1}
\dots s_{n_{l-2}}} + \dots + k_d^{n_l} \underline{s_{n_1 -1} \dots
s_0}$$

\noindent where $\underline{s_a \dots s_b}$ is the number with digits
$s_a, \dots, s_b$ respectively (when written from left to right, in base
$B$).\\

\begin{cor}\label{C1}
$\aaa_{n_1, \dots, n_l,s}(k_d) \equiv B^{-n_l}s \mod d$.
\end{cor}

\begin{proof}
We compute, for general $i$,
$$k_d^{n_l - n_i} \underline{s_{n_{i+1}-1} \dots s_{n_i}} = k_d^{n_l -
n_i} \sum_{j=n_i}^{n_{i+1}-1} s_j B^{j - n_i} \equiv \sum_{j =
n_i}^{n_{i+1} - 1} s_j k_d^{n_l - n_i - j + n_i} \equiv \sum_{j =
n_i}^{n_{i+1} - 1} s_j k_d^{n_l - j}$$
Therefore we have
$$\aaa_{n_1, \dots, n_l,s}(k_d) \equiv \ssumm k_d^{n_l-j} \equiv
k_d^{n_l} \ssumm k_d^{-j} \equiv B^{-n_l} \ssumm B^j \equiv B^{-n_l}s$$
(where everything is modulo $d$). Hence we are done.
\end{proof}\hfill

\noindent {\bf Remarks:} Thus, the number $k_d$ (as above) also satisfies
(for any $s$) :
\begin{enumerate}
\item For all such tuples $(n_1, \dots, n_l),\ s \equiv 0 \mod d \text{
iff } \aaa_{n_1,\dots,n_l,s}(k_d) \equiv 0 \mod d$.
\item Hence for all $i,\ s \equiv 0 \mod d \text{ iff } \aaa_{i,s}(k_d) =
\underline{s_m \dots s_i} + k_d^i \underline{s_{i-1} \dots s_0} \equiv 0
\mod d$.
\item Setting $n_i = i\ \forall i= 1, \dots, m-1$ in (1), we get that $s
\equiv 0 \mod d\\ \text{ iff } \aaa_{1,2, \dots, m-1, s}(k_d) =
\aaa_s(k_d) \equiv 0 \mod d$.\\
\end{enumerate}
In particular, the ``reverse" divisibility tests (in Theorem \ref{Tprep1}
and Section \ref{secKh}) mentioned above, do hold.\\

\noindent {\bf Example:} There is a well-known divisibility test for 7
that involves splitting up numbers into groups of three digits, and then
taking the alternate sum and difference of the numbers. For instance, to
check 142857142 for divisibility by 7, we compute: 142 - 857 + 142 = 284
- 857 = -573 $\equiv -13 \equiv 1 \mod 7$, whence $7 \nmid 142857142$.
The reason this works is the corollary above (and the remarks following
it), if we note that 1000 $\equiv -1 \mod 7$, so that $10^6 \equiv 1 \mod
7$.\\

We now show that the $k$-values of terms in an arithmetic progression
(more precisely, having the same ``unit's digit"), themselves form an
arithmetic progression. This is useful to calculate $k$-values for
general $d$, if we know the $k$-value of only the unit's digit of $d$.
(For example, knowing $k_9$, we can easily calculate $k_{19}, k_{29}$,
etc.).\\

\begin{theorem}\label{Tprep2}
Take $a \in E$ ($a,B$ coprime), and $m \in E$. Let $d=Bm+a$. Then
\begin{equation}\label{E6}
k_d = k_{a+mB} \equiv k_a + ml \mod d, \text{ where } l = \frac{1}{a}
(Bk_a - 1)
\end{equation}
\end{theorem}

\begin{proof}
We compute: $B(k_a+ml) = Bk_a + Bml = (1+al) + Bml = 1 + l(a+mB) \equiv
1 \mod (a+mB) \equiv 1 \mod d$. Hence by definition of $k_d$, we are
done.
\end{proof}\hfill

\begin{lemma}\label{Lprep3}
Let $d | (B^t - a)$ where $a \in E^\times$ and $t \in \nn$. Then $k_d
\equiv B^t(aB)^{-1} \mod d$.
\end{lemma}

\begin{proof}
$B (B^t(aB)^{-1}) = a^{-1}(B^t) = a^{-1}(B^t - a) + 1 \equiv 1 \mod d
\equiv B k_d \mod d$. Hence we cancel $B$ (modulo $d$) to get the
result.
\end{proof}

This (last) result tells us, in particular, the divisibility tests for 3,
9, and 11. We now give two examples in the Euclidean domain $E=F[X]$
(where $F$ is a field.)\\

\noindent {\bf Example 1:} 
Let us work with $B = X,\ d = X^2+1$. We want to check whether or not $s
= X^4 + X^2$ is divisible by $d$. Clearly, $s = 10100$ in base $X$, and
to apply our test, we note that $(-X)d = -X^2 \equiv 1 \mod d$. Hence $k
= -X = -10$. So we get $\alpha_s(k) = 1 + (-10)^2 = 101$ and this is
indeed divisible by $101 = d$. Therefore, so is the original $s$.\\

\noindent {\bf Example 2:}
We now show how the Factor Theorem for polynomials is an example of the
reverse test.
We still work with $B = X$, and $d = X-c$ for some $c \in F^\times$. We
want to check if a general element
\[s = \underline{s_m s_{m-1} \dots s_1 s_0} = \ssumm X^j\]

\noindent is divisible by $d$ or not. From Lemma \ref{Lprep3} above, $k_d
\equiv c^{-1} \mod d$. Thus to check whether or not $d|s$ is the same as
checking whether or not $(X-c)$ divides $\alpha_s(k_d) = \ssumm
(c^{-1})^{m-j} \in F$. Since this is a scalar, it should be zero, and
multiplying by $c^m$ should still give zero. Thus,
\[ (X-c) | s(X) \mbox{ iff } \ssumm c^j = 0 \mbox{ iff } s(c) = 0 \]

\noindent and this is just the Factor Theorem for polynomials. Moreover,
Theorem \ref{Tprep2} above says that the ``$k$-value" function (taking
$d_m = mX-c$ to $k(d_m)$), is a linear function in $m \in F$.\\

\subsection{Testing for divisibility by general $n \in E$}\hfill\break

We now check whether or not a given number $s$ is divisible by $n$, where
$n \in E$ need not be coprime to $B$, the base in which both $s$ and $n$
are written. To do this, we take $n = d p_1^{\aaa_1} p_2^{\aaa_2} \dots
p_h^{\aaa_h}$ (by unique factorization in the UFD $E$), where $d$ is
coprime to $B$, and $p_1, p_2, \dots, p_h$ are some or all of the primes
that divide $B\ (\aaa_i \in \N \text{ for } i=1,2,\dots,h)$.

Divisibility by $d$ is tested by the method given above, and for each
$p_i$, we perform the analogue (in base $B$) of the ``$n$-right-digits
test":\\

If $\beta_i$ is the highest power of $p_i$ dividing $B$, then to test for
$p_i^{\aaa_i}$, we check whether or not the number formed by the
$\grint{\aaa_i / \beta_i}$ rightmost digits of $s$ ({\it in that order}
and in base $B$) is divisible by $p_i^{\aaa_i}$. (Here, $\grint{x}$ is
the least integer greater than or equal to $x$, e.g. if $\beta_i |
\aaa_i$, then $\grint{\aaa_i / \beta_i} = \aaa_i / \beta_i$.) This is
because the place value of each of the other digits is $\geq
B^{\grint{\aaa_i / \beta_i}}$, and this is already divisible by
$p_i^{\beta_i \grint{\aaa_i / \beta_i}}$, which is a multiple of
$p_i^{\aaa_i}$.

For example, take $E=\Z$. In base 10, to check for divisibility by 125,
we only need to check the three rightmost digits - they have to be 000,
125, 250, ..., 750, or 875. As another example, consider $B = 12$. Since
$12 = 3 \times 2^2$, hence to check divisibility by $8=2^3$, the
$\grint{3/2} = \grint{1.5} = 2$ rightmost digits have to be checked.\\

\section*{\bf Part 2 : Recurring Decimal Representations}

\section{Decimal representation of fractions - a few observations}

We know that every rational number can be represented as a number having
a terminating or a recurring set of digits in its decimal representation.
We observe certain interesting properties of the digits in the recurring
portion, e.g. $\frac{1}{7} = 0.142857\ 142857\ \dots,$ and upon
multiplying 142857 by any number from 1 to 6, we get the same six digits
in the same cyclic order (142857, 285714, 428571, 571428, 714285,
857142). On multiplying it by 7, we get 999999, so that in the decimal
representation, it becomes $0.9999\dots = 1$. Another similar case is
$\frac{1}{17} = 0.05882352\ 94117647\ 0588\dots$, and if $0588\dots7647$
(including the zero) is multiplied by any number from 1 to 16, then the
same sixteen digits are obtained in the same cyclic order; if it is
multiplied by 17 then the result, of course, is $9999\dots9999$ (sixteen
9's). This is further addressed in the comments following Lemma
\ref{Lsquare} below.\\

We also observe the following phenomenon, which we quote from \cite{Kv}:
We know that $1/7 = 0.\overline{142857}$. If we square 142857, we get
20408122449. Take the first six digits (from the right) and add them to
the number formed by the rest (we take six digits because 142857 has six
digits). We then get 20408+122449 = 142857.

Similarly, if we take $1/27 = 0.\overline{037}$ and square 37 then we get
1369, and adding the rest to the three rightmost digits yields 1+369 =
370, again a multiple of the original number. We shall also show this
phenomenon in general, in Lemma \ref{Lsquare} below.\\

We observe another interesting property in both cases - and more
generally, in the case of every prime whose reciprocal has an even number
of digits in its recurring decimal portion. Namely, if the second half of
numbers is kept below the first half, then the sum of every pair of
corresponding digits is 9 (e.g. in 142 857, 1+8 = 4+5 = 2 + 7 = 9). We
give further examples below, and explain this in Lemma \ref{Leven} and
Theorem \ref{Tnorec} below.\\

First of all, how do we get recurring decimals ? For an odd prime $p \neq
5$, we look at the {\it order} $b$ of $10$ in $(\Z / p \Z)^\times$. We
have $10^b \equiv 1 \mod p$. Thus $p$ divides $10^b-1 = 999 \dots 9$ ($b$
nines). Let $10^b-1$ be denoted by $9_b$, and let $\frac{9_b}{p} = r_p$.
Thus, we have $\frac{1}{p} = \frac{1}{p}(0.999\dots) = \frac{1}{p}(0.9_b
9_b \dots) = 0.r_p r_p \dots$. Then the number of digits in the recurring
part of the reciprocal of $p$ gives the number of terms (i.e. $b$) in the
chain. Related facts are shown in Lemma \ref{Lbasic} below.\\

Given $p$, we have various {\it chains} of recurring decimals - for
instance, we have the (unique) chain $\overline{142857} =
\overline{571428}$ etc. for 1/7. However, for some numbers (e.g. 13), we
might have more than one chain. Let $c_p$ be the number of distinct
chains. Then $b_p c_p = p-1$.

Now consider the situation in base $B$ (for any $B>1$). In what follows,
the phrase ``decimal expansion" also refers to the case of a general base
$B$. In general, for composite $d$, the recurring parts of $t/d$ might
have recurring chains of different lengths as well (e.g. for $d=21$, in
base 10, the chains for 7/21 and 14/21 have length 1 each, and are
different from the 6-digit chains for the other residues). However, if we
denote the different lengths of all possible chains by $b_i$, and the
number of chains of a given length ($=b_i$) by $c_i$ (for $i = 1, 2,
\dots$), then we have
\begin{equation}\label{eqnbc}
\text{ For any } d\ (d \in \N,\ d \text{ coprime to } B),\ \sum_i b_i
c_i = d
\end{equation}
We also remark that the number of chains $c_p$ is also known as the {\it
residual index}, and it clearly equals $[(\Z/p\Z)^\times : \langle B
\rangle]$.

Note that in base 10, we take $d/d = 1 = 0.999 \dots = 0.\bar{9}$, and
this makes the above sum one more than the $d-1$, that we obtained above,
in the case when $d$ was a prime. Similarly for base $B$.

Of course, we can find many composite numbers $d$ which have only one
value of $b$ ($= b_d$, say), e.g. $b_{91} = 6$. This happens, for
instance, when $d$ is the product of two or more primes, all distinct,
and all having the same $b$-value. This is addressed in Proposition
\ref{Pallsame} below.\\

\begin{proof}[Proof of Equation \eqref{eqnbc}]
Look at $\{ t/d : 1 \leq t \leq d \}$, and look at their ``dec"imal
representations. Clearly, there are no repetitions among the various
decimals, for the numbers in the above set are all distinct (modulo $d$).
Further, the numbers in the set whose recurring chains coincide with the
chain $t/d = 0.\overline{a_1 \dots a_b}$ are precisely $t/d,\ (Bt \mod
d)/d,\ (B^2t \mod d)/d, \dots,\ (B^{b-1}t \mod d)/d$. Hence, given any
chain, {\it every} cyclic permutation of it corresponds to $t/d$ for some
$d$.

This gives a bijection between all cyclic permutations of a fixed chain
($b_i$ of them), and $b_i$ of the residue classes modulo $d$. Hence we
get equation \eqref{eqnbc} above, because there are $d$ residue classes
modulo $d$ (where we also include $1 = d/d = 0.\overline{9}$ etc.).
\end{proof}\hfill

\section{Properties of recurring decimals - proofs for a Euclidean
domain}

We now work on fractions and recurring ``dec"imals in base $B$, in a
Euclidean domain $(E,\nu)$. Though we cannot talk about convergence here,
we can still talk about repeating decimals representing a fraction, and
such things.

Observe that if $0 \neq d$ is not a unit, then $\nu(1) < \nu(d)$, so
$1/d$ should have an expansion of the form $0.a_1 a_2 \dots = \sum_i a_i
B^{-i} \in k(E)$, the quotient field of $E$. (Actually, we work in
$E[[B^{-1}]] = E[[X]] / (BX-1)$.) We now find out what the $a_i$'s are.\\

\noindent {\bf Remarks :}
\begin{enumerate}
\item We do not worry about convergence issues here, and only associate
to a given fraction, a {\it sequence} $a_i$ of numbers, without worrying
whether these numbers are actually digits or not (i.e. whether $\nu(a_i)
\geq$ or $< \nu(B)$). Note that if $E=\Z$, then all sequences are
actually decimal expansions that converge etc.\\

\item If $d \neq 0$ and $d \nmid a$ in $E$, then $a/d = 0.a_1 a_2 \dots$
iff $\nu(a) < \nu(d)$. In other words, there exists a representation of
$a/d$ as $0.a_1 a_2 \dots$ iff the quotient (when we divide $a$ by $d$ in
$E$) is zero, iff $a = 0 \cdot d + r = r$ satisfies the Euclidean
algorithm, namely, that $\nu(a) = \nu(r) < \nu(d)$.\\
\end{enumerate}

\begin{lemma}\label{Lprep4}
Suppose $0 \neq d$ is not a unit, and we define $\frac{1}{d} = 0.a_1 a_2
\dots$. Then $a_i = \leint{B^i/d} - B \leint{B^{i-1}/d}$.
\end{lemma}

\noindent Here, by $\leint{a/b}$, we mean the quotient when we divide $a$
by $b$.\\

\noindent {\bf Warning} : We need to choose and {\it fix} the quotients
$\leint{B^i/d}$, since these are not, in general, unique ! Moreover,
convergence cannot be discussed, and the $a_i$'s {\it need not be digits}
(i.e. we do not know if we have $\nu(a_i) < \nu(B)$ or not).

\begin{proof}
Suppose we write $1/d = 0.a_1 a_2 \dots$. Clearly, the quotient (i.e. the
number to the left of the decimal point) is zero because of the Euclidean
algorithm (and since $\nu(1) < \nu(d)$). Next, $a_1$ is clearly of the
desired form, because of the same reason.

We now show that the $a_i$'s are as claimed, by induction on $i$. The
base case was done above. To show the claim for $a_i$, multiply the
equation above by $B^i$. Thus, $B^i/d = a_1 a_2 \dots a_i.a_{i+1}a_{i+2}
\dots$, so taking the quotient of both sides yields
$$\leint{B^i/d} = \underline{a_1 a_2 \dots a_i} = \sum_{j=1}^i B^{i-j}
a_j$$

This is because $B^i = d(\underline{a_1 \dots a_i}) +
d(0.\overline{a_{i+1} \dots a_n a_1 \dots a_i}) = dq+r$, say (where $r =
d(0.\overline{a_{i+1} \dots a_n a_1 \dots a_i}) \in E$). By the above
assumption (since $(B,d) = 1$), we get that $\nu(r) < \nu(d)$, so $q$ and
$r$ are clearly the quotient and remainder, as desired.\\

Assume by induction that we know the results for $a_1, \dots, a_{i-1}$.
Then
$$a_i = \leint{B^i/d} - \sum_{j=1}^{i-1} B^{i-j}(\leint{B^j/d} -
B\leint{B^{j-1}/d})$$

and the latter is a telescoping sum, so we get
$$a_i = \leint{B^i/d} - B \leint{B^{i-1}/d} + B^i \leint{B^0/d}$$

But the last term is zero, because $\nu(B^0) < \nu(d)$. Hence we are
done.
\end{proof}\hfill

\begin{lemma}\label{Lvalue}
The recurring decimal $0. \overline{a_1 \dots a_n}$ denotes the fraction
$\underline{a_1 \dots a_n} / (B^n-1)$. If $\underline{a_1 \dots a_n} =
B^n-1$, then the decimal equals 1.
\end{lemma}

\begin{proof}
Suppose $e = 0. \overline {a_1 \dots a_n}$. Then $B^n e = \underline{a_1
\dots a_n}.\overline{a_1 \dots a_n}$. Subtracting from this the
definition of $e$, we are done. If, moreover, $\underline{a_1 \dots a_n}
= B^n - 1$, then $(B^n - 1)(e-1) = 0$. Since $B$ is not a unit, $B^n \neq
1$; since $E$ is an integral domain, we conclude that $e=1$.
\end{proof}\hfill

Next, to talk about periodicity of any and every $a/d$ (for fixed $d$
coprime to $B$) in such a decimal expansion, we need the following

\noindent {\bf Standing Assumption} : $B$ has finite order in
$(E/dE)^\times$. (In particular, $d \notin E^\times$.)\\

\noindent Denote this order by $e$. Thus $e = o_d(B) = o_d(B^{-1} \mod d)
= o_d(k_d)$. We now get

$$\frac{1}{d} = \frac{B^e-1}{d} \bigl( \frac{1}{B^e} + \frac{1}{B^{2e}} +
\dots \bigr) = 0.r_d r_d \dots, \mbox{ where } r_d = \frac{B^e-1}{d}$$

\noindent and this is how we get recurring decimal expansions. Moreover,
from the proof of Lemma \ref{Lprep4} above, $r_d = (B^e-1)/d =
\leint{B^e/d}$ (because $\nu(1) < \nu(d)$) $ = \underline{a_1 \dots
a_e}$, and hence $r_d$ has $e$ digits. Thus, the correct way to look at
this, in a Euclidean domain, is to look at the {\it order} of $B$ in
$(E/dE)^\times$, instead of the recurring decimal expansion or the number
of digits therein.

\begin{lemma}\label{Lbasic}
If $t$ is coprime to $d$, and $B$ has finite order $e$ in $(E/dE)^\times$,
then the sequence of $a_i$'s associated to $t/d$ is recurring with period
$e$.
\end{lemma}

\begin{proof}
The length of the sequence is the smallest $e' \in \N$ such that $t
B^{e'} \equiv t \mod d$. Since $t$ is coprime to $d$, we get $e'=e$ by
definition.
\end{proof}\hfill

\noindent We now turn to repeating decimals with even period. We have

\begin{lemma}\label{Leveniff}
$B$ has (finite, and) even period, modulo a prime $p \nmid 2B$, iff $p |
(B^l+1)$ for some $l \in \N$.
\end{lemma}

\begin{proof}
Suppose $p | (B^l+1)$ for some positive $l$; we assume, moreover, that
$l$ is the least such. Thus $p | (B^{2l} - 1)$, whence the order $o_p(B)$
is finite and divides $2l$. Since $p \nmid 2$, hence $p \nmid B^l-1$,
whence $o_p(B) \nmid l$. Thus $o_p(B) = 2m$ for some $m | l$.

Conversely, if $o_p(B) = 2m$, then $p | (B^m+1)(B^m-1)$. By definition of
order, $p \nmid (B^m-1)$, whence $p | (B^m+1)$, and we are done. (Note,
moreover, that continuing the above proof of the first part, we find that
$l=m$ by choice of $l$.)
\end{proof}\hfill

\begin{lemma}\label{Leven}
Suppose $d | (B^l+1)$ for some $l \in \N$. Then for all $a$ coprime to
$d$, the chains starting from the $i$th and the $(l+i)$th ``digits" in
the recurring part of $a / d$, add up to a recurring chain that
represents 1 or 0.
\end{lemma}

\begin{proof}
Since $B^l \equiv -1 \mod d$, hence the chain corresponding to adding up
the fractions $(aB^i \mod d)/d$ and $(aB^{l+i} \mod d)/d$, yields the
fraction $((aB^i + a(-1)B^i) \mod d) / d = d/d = 1$, as claimed.
\end{proof}

\noindent {\bf Examples}:
\begin{enumerate}
\item Suppose $E=\Z$ and $B=10$, so that $\nu(B) > \nu(B-1)$. Then $1/7 =
0.\overline{142\ 857}$, which satisfies the above assumptions, since the
chains add up to $0.\overline{999999} = 0.\overline{9} = 1$.

Similarly, $1/11 = 0.\overline{09}$. Here's another way to write this.
Observe that $9 = 10 \cdot 0 + 9 = 10 \cdot 1 + (-1)$, with $\nu(-1) =
|-1| < \nu(10)$. Hence we can also write $1/11 = 0.\overline{1(-1)}$, and
then $1 + (-1) = 0$. This corresponds to the recurring chain
$1.\overline{0}$ for 1.

\item As another example, consider $B = F[X]$, with char $F \neq 2$.
Suppose $B=X$ and $d = X^2+1$ (so that $\nu(B) = \nu(B-1)$). Clearly, $B$
has order 4 modulo $d$, since $d | (X^4-1)$. So $1/d$ should have four
digits. But we also know what these are : $r_d = (B^e-1)/d = (X^4-1) /
(X^2+1) = X^2 - 1 = 010(-1)$ in base $B=X$. Thus $1/d =
0.\overline{010(-1)}$ and once again, every pair of alternate digits adds
up to zero.
\end{enumerate}

\noindent We thus see that the case $E=\Z$ and $1 = 0.\overline{9}$ is
somewhat of an ``exception" ! This is made more precise now.\\

\begin{theorem}\label{Tnorec}
Given a Euclidean domain $(E,\nu)$ that is not a field, the following are
equivalent:
\begin{enumerate}
\item $\nu(B) = \nu(B-1)\ \forall B \notin E_0$.

\item $\nu(B) = \nu(B-c)\ \forall B \notin E_0, c \in E_0$.

\item The Euclidean algorithm yields a unique quotient and remainder when
we divide any $B^n-c$ by $B$, for any $n \in \N, B \notin E_0, c \in
E_0$.

\item There is no $B \notin E_0$, so that 1 can be expressed (in base
$B$) as a recurring decimal $0. \overline{a_1 \dots a_n}$ for some $n \in
\N$, with $\nu(a_i) < \nu(B)$ or $a_i = 0$.
\end{enumerate}
\end{theorem}

\begin{proof}
We prove a series of cyclic implications.\\
\noindent{$(1) \Rightarrow (2)$:}\\
This is clear for $c=0$, and for $c \in E^\times$, we have (assuming (1)
and using part (2) of Lemma \ref{L1})
\[ \nu(B-c) = \nu(c^{-1}(B-c)) = \nu(c^{-1}B - 1) = \nu(c^{-1}B) = \nu(B)
\]

\noindent{$(2) \Rightarrow (3)$:}\\
Suppose $B^n - c = qB + d$. Then $B | (c+d)$, whence $(c+d)$ is not a
unit. If $d \neq -c$ then $c+d$ is a nonzero nonunit (whence $d \neq 0$).
Hence $\nu(c+d) = \nu(d)$ by (2). Moreover, $B | (c+d)$, so $\nu(B) \leq
\nu(c+d) = \nu(d) < \nu(B)$, the last inequality a consequence of the
Euclidean algorithm. This gives a contradiction, whence $d = -c$ and
hence $q = B^{n-1}$, as desired.\\

\noindent{$(3) \Rightarrow (4)$:}\\
Suppose we can write 1 as a recurring ``decimal". Then by Lemma
\ref{Lvalue}, $B^n-1 = \underline{a_1 \dots a_n} = a_n + B \cdot
\underline{a_1 \dots a_{n-1}}$. Assuming (3), one notes that $a_n = -1$,
since $-1 \in E_0$. Thus $\underline{a_1 \dots a_{n-1}} = B^{n-1}$. Once
again, by (3) we get that $a_{n-1} = 0$. Proceeding inductively, we get
that $a_i = 0\ \forall 1 < i < n$, and finally, we are left with $a_1 =
B$, a contradiction to $\nu(a_1) < \nu(B)$. Hence 1 cannot be written in
recurring form.\\

\noindent{$(4) \Rightarrow (1)$:}\\
We prove the contrapositive. Suppose we have a $B \notin E_0$ so that
$\nu(B) \neq \nu(B-1)$. There are two cases. If $\nu(B) > \nu(B-1)$, then
(using Lemma \ref{Lvalue}) we have the expansion $1 = 0.\overline{(B-1)}$
in base $B$. On the other hand, if $\nu(B) < \nu(B-1)$, then the
expansion $e = 0.\overline{B (-B)}$ in base $(B-1)$, again equals 1 (by
Lemma \ref{Lvalue}).
\end{proof}\hfill

\begin{lemma}\label{Lsquare}
Suppose $B^e -1 = d \cdot r_d$, where $r_d$ is the recurring part. Fix
$k,l \in \N$ and consider $k r_d$. Then the number formed by the first
$le$ ``digits" (or $a_i$'s), added to the number formed by the rest,
yields a number that is divisible by $r_d$.
\end{lemma}

\begin{proof}
The number in question is obtained as follows: suppose the ``rest of the
digits" form a number $s'$. (For instance, if $E=\Z$, then we could write
$s'=\leint{k r_d / B^{le}}$, where this denotes the greatest integer in the
quotient.) Then the total number obtained is
\[ s' + (k r_d - B^{le} s' )= k r_d - (B^{le}-1) s' \]

\noindent and this is divisible by $r_d$, the quotient being $k - ds'(1 +
B^e + \dots + B^{le-e})$.
\end{proof}

We now show some of the facts from the previous sections. First of all,
observe that if $t/d = 0.a_1 \dots a_i \dots$, then $0.a_{i+1} a_{i+2}
\dots$ comes from the expansion of $B^i t / d$, or from $(B^i t \mod d) /
d = a/d$, say. Thus if $r_d = (B^e - 1)/d$, then the recurring part of $a
/ d$ would correspond to $a r_d$. This explains why recurring decimal
expansions keep getting permuted cyclically, upon being multiplied by
various numbers.\\

Next, we define $D_d$, for any $d$ satisfying the above assumption, to be
the order $e$ of $B$ in $(E / dE)^\times$, or equivalently, the ``number
of terms" of the sequence $(a_i)$, in $r_d$ above.\\

Suppose $d = \prod_{i=1}^n p_i^{\alpha_i}$ as above, where $p_i$ are
primes in $E$ (which is a unique factorization domain). If $d'|d$, then
$$0 \to dE \to d'E\ (\hookrightarrow E)$$

\noindent and hence
$$E/dE \to E/d'E \to 0$$

\noindent (The kernel is, of course, $d'E/dE$, so we get the short exact
sequence
$$0 \to d'E/dE \to E/dE \to E/d'E \to 0$$

\noindent for every choice of $d'$.) In particular, $B$ has finite order
in every $(E/d'E)^\times$ (because if $B^e \equiv 1 \mod d$, then $B^e
\equiv 1 \mod d'$); in particular, in $(E/p_iE)^\times$ as well. Thus we
can talk of $D_{p_i}$ as well. However, this does not guarantee that the
standing assumption holds for $p_i^n$  for all $n$. So we need to make a
stronger assumption.\\

However, let us remark that once $B$ has finite order in $(E/dE)^\times$,
every $t/d$ has a recurring decimal expansion. For, we write $t/d =
t'/d'$ in lowest terms, and then from the above remarks, $B$ has finite
order in $(E/d'E)^\times$, and we are done.\\

\begin{prop}\label{Pallsame}
Suppose $d = p_1 \dots p_h$ for distinct primes $p_i$, all coprime to
$B$. Let $b_i\ (\in \N)$ denote the order of $B$ in $(E/p_iE)^\times$.
\begin{enumerate}
\item Then $B$ has finite order in $(E/dE)^\times$, and this order,
denoted by $b_d$, say, is the l.c.m. of the $b_i$'s.
\item If $b_i = b_d\ \forall i$, then $B$ has order $b_d$ in
$(E/d'E)^\times$, for {\emph every} nonunit $d' | d$.
\end{enumerate}
\end{prop}

\begin{proof}
\hfill
\begin{enumerate}
\item Clearly, $B^{lcm} -1 \equiv 0 \mod p_i\ \forall i$, where $lcm$
denotes the l.c.m. of the $b_i$'s. Since the $p_i$'s are mutually
coprime, and $E$ is a UFD, hence $\prod_i p_i$ also divides $B^{lcm} -
1$. Thus $B$ has finite order in $(E/dE)^\times$, and moreover, this
order has to divide the l.c.m. But the l.c.m. has to divide the order as
well, so we are done.

\item Suppose $b_i = b_d\ \forall i$. Now, if $d \nmid a$, then $a/d$ is
of the form $r/s$, where $r$ and $s$ are coprime, and $s$ is a product of
some of the $p_i$'s. Clearly, $r/s$ and $1/s$ have the same number of
recurring digits, because $(r,s) = 1$. But $s \neq 1$ is a product of
some $p_i$'s, and each $1 / p_i$ has $b_d$ recurring digits. Hence so
does $1/s$, and so does $r/s = a/d$.
\end{enumerate}
\end{proof}

\section{The period of recurrence of a fraction in a Euclidean
domain}\label{sec9}

\noindent We have the following

\begin{lemma}\label{Lprep9}
If $B$ has finite order $e$ modulo $p$, then the following are
equivalent:
\begin{enumerate}
\item $B$ has finite order modulo $p^n$ for any $n$.
\item There exists a unique \emph{prime integer} $p' \in \Z$ (which is
prime in $\Z$, but may not be prime in $E$), so that $p|p'$.
\item There exists $n \in \Z$ so that $p | n$.
\end{enumerate}
\end{lemma}

\noindent {\bf Remark}: By $\Z$ we mean the image of $\Z$ in $E$, for we
do not know if $E$ has characteristic zero or not.\\

\begin{proof}
Suppose $B$ has order $e$ modulo $p$. Thus $B^e = 1 + \alpha p$ for some
$\alpha \in E$. Suppose $\alpha = p^f \beta$, where $p \nmid \beta$. Now,
$B^{en} \equiv 1 + n \alpha p \mod p^{f+2}$ (ignoring higher order terms).
Thus, $B$ has finite order modulo $p^{f+2}$ iff there is some $n$ such that
$p^{f+2} | n \alpha p$. But from our assumptions on $f$, this means $p |
n$. Thus $(1) \Rightarrow (3)$.

Conversely, assume that $p | n$. Now observe that if $s \equiv 1 \mod
p^m$ for any $s \in E,\ m \in \N$, then $s = 1 + \beta p^m$ for some
$\beta \in E$, whence
\[ s^n = (1 + \beta p^m)^n \equiv 1 + \beta n p^m \mod p^{2m} \equiv 1 +
\beta n p^m \mod p^{m+1} \equiv 1 \mod p^{m+1} \]

\noindent since $p | n$ and $m \in \N$. In particular, since $B^e \equiv
1 \mod p^{f+1} \equiv 1 \mod p$, hence $B^{e n^{m-1}} \equiv 1 \mod p^m$
for all $m \in \N$. Hence $(3) \Rightarrow (1)$.

That (2) and (3) are equivalent is easy : in one direction, choose
$n=p'$. In the other direction, decompose $n$ into prime factors (as {\it
integers}), and note that in $E$, $p$ must divide some integer prime
factor of $n$. If $p$ divided two such, then it would divide any integer
linear combination of these prime numbers, whence $p$ would divide 1, and
hence be a unit. This is not possible.
\end{proof}

\noindent Thus, we now make the following\\

\noindent {\bf Standing Assumptions} :
\begin{enumerate}
\item $(E,\nu)$ is a Euclidean domain, not a field. Fix $B \notin E_0$.

\item For each prime $p \in E,\ p \nmid B$, (a) $B$ has finite order
modulo $p$, and (b) there exists a unique {\it prime integer} $p' \in \Z$
so that $p | p'$ (and $p' \neq 0$) in $E$.
\end{enumerate}

\noindent As an easy consequence, we have

\begin{lemma}
$E$ (or its quotient field) has characteristic 0.
\end{lemma}

\begin{proof}
$E$ is not a field, so there exists a nonzero nonunit in $E$. By unique
factorization, there exists a nonzero prime $p \in E$, so by assumption
there exists a nonzero prime integer $p'$ that is a multiple of $p$. Now,
$p' \notin E_0$ since $p \notin E_0$. This would give a contradiction if
$E$ had positive characteristic, so we are done.
\end{proof}\hfill

\begin{lemma}\label{Luniqueint}
If $(B,d) = 1$, then $B$ has finite order modulo $d$. Further, there is
a unique $n \in \N$ such that for each $m \in \Z$, $d | m$ in $E$ iff $n
| m$ in $\Z$.
\end{lemma}

\begin{proof}
We argue by induction on the number of distinct prime factors of $d$
(ignoring multiplicities). If this number is 1, then $d$ is a prime
power, and we are done by the previous lemma.

Suppose we now have a general $d = p^r d'$ for some $p$ prime in $E$,
coprime to $d'$, and suppose we know the result for $d'$. Then $B^f = 1 +
\alpha d'$ for some $\alpha, f$. If $B^e \equiv 1 \mod p^r$, then $B^{ef}
\equiv 1$ modulo both $d'$ and $p^r$, which are coprime. Hence by unique
factorization, $B^{ef} \equiv 1 \mod d$, and we are done.

Finally, if $d = \prod_i p_i^{\alpha_i}$, and $p_i | p_i' \in \Z$, then
$d | \prod_i (p'_i)^{\alpha_i} \in \Z$. Thus the set $S$ of integers
divisible by $d$, properly contains the element 0 (since $E$ has
characteristic zero, from above). Moreover, $S$ is clearly an ideal in
$\Z$, hence is principal and generated by some $n \in \N$. This is the
required $n$ (and $d|n$), and we are done.
\end{proof}\hfill

We now calculate $D_d$ = the number of digits in the recurring decimal
part of $\frac{x}{d}$, where $d$ is coprime to both $x$ and to $B$, and
$d = p_1^{\aaa_1} p_2^{\aaa_2} \dots p_h^{\aaa_h}$ (the $p_i$'s here are
distinct primes coprime to $B$, and $\aaa_i \in \N$).

For the rest of this section, fix $d = \prod_i p_i^{\alpha_i}$ coprime to
$B$, and fix $p$, a prime factor of $d$. Define $q_i(p) = q_i = D_{p^i}$
= order of $B$ in $(E/p^iE)^\times$, and define $f$ to be the highest
power of $p$ that divides the {\it prime integer} $p'$ (as in the
Standing Assumptions above). Also let $\nu_p(m)$ be the largest power of
$p$ that divides $m$ in the UFD $E$.\\

We first prove the following generalization of \cite[P.1.2 (iv),
pp.11]{Ri} (special cases of the result in \cite{Ri} were known to
Euler):

\begin{prop}\label{Peuler}
Given $x \neq y$ in $E$, define $a(m) = \fract{x}{y}{m}$ for $m \geq 0$.
Now assume that $p | (x-y)$ and $p \nmid y$. Then $p | a(m)
\Leftrightarrow p' | m$. If, moreover, $\nu_p(x-y) > \leint{f/(p'-1)} \in
\N_0$, then for all $m \in \N$, we have
\[ \nu_p(a(m)) = \nu_p \bigg( \fract{x}{y}{m} \bigg) = \nu_p(m) \]

\noindent In particular, $(x^{p'}-y^{p'})/(x-y)$ is divisible by $p^f$
but not by $p^{f+1}$.\\
\end{prop}

\begin{proof}
Firstly, it is easy to check that $(x-y) | (a(m) - my^{m-1})$ for all
$m$. Now suppose $p | (x-y)$. Then $a(m) \equiv my^{m-1} \mod (x-y)$,
whence the same holds modulo $p$. Since $p,y$ are coprime, we conclude
that $p | a(m) \Leftrightarrow p | m \Leftrightarrow p' | m$, by Lemma
\ref{Luniqueint} above.

Now suppose that $\nu_p(x-y) > f / (p'-1)$. Setting $x-y=t$, we get that
$a(p') = ((y+t)^{p'} - y^{p'}) / t = \sum_{i=1}^{p'} \binom{p'}{i}
y^{p'-i} t^{i-1}$. Over here, every term except the first and last one,
is divisible by $p' \cdot p$, since there is a binomial coefficient and a
power of $t$ in each term. Thus, modulo $p^{f+1}$, we have $a(p') \equiv
p' y^{p'-1} + t^{p'-1}$. The last term is divisible by $p^{n(p'-1)}$,
where $n = \nu_p(x-y) > f/(p'-1)$; hence it is divisible by $p^{f+1}$.
The first term is only divisible by $p^f$, since $p,y$ are coprime. Thus
the last line of the result is proved.\\

To prove the rest of the result, suppose $m = u {p'}^s$, where $p' \nmid
u$. Now define $X_i = x^{(p')^i}, Y_i = y^{(p')^i}$, and note that
\[ \fract{x}{y}{m} = \fract{X_s}{Y_s}{u} \cdot \prod_{i=0}^{s-1}
\frac{X_{i+1} - Y_{i+1}}{X_i - Y_i} = \fract{X_s}{Y_s}{u} \cdot
\prod_{i=0}^{s-1} \fract{X_i}{Y_i}{p'} \]

\noindent where the product on the right telescopes. Since $(x-y) | (X_i
- Y_i)$, hence $\nu_p(X_i - Y_i) > f / (p'-1)$ for all $i$. Hence by the
previous paragraph, each term in the product is exactly divisible by
$p^f$; moreover, the remaining term is not divisible by $p$ by the first
part of this result, since $p \nmid u$. Hence $\nu_p(\fract{x}{y}{m}) =
fs = \nu_p(m)$.
\end{proof}\hfill

\begin{theorem}
For all $n \in \N,\ q_{n+1} = q_n$ or $p'q_n$. Now let $g$ be the least
natural number so that (a) $g > f / (p'-1)$, and (b) $q_{g+1} = p'q_g$.
Then for all $n \geq 0$, we have
\[ q_{g+n} = (p')^{\grint{n/f}} q_g \]
\end{theorem}

\begin{proof}
First, the general case. We assume $q_{n+1} \neq q_n$. Now, $B^{q_{n+1}}
\equiv 1 \mod p^{n+1} \equiv 1 \mod p^n$, whence $q_n | q_{n+1}$. So let
$q_{n+1} = mq_n$. We now use the first part of Proposition \ref{Peuler},
setting $x = B^{q_n}$ and $y=1$. Thus, the order of $B^{q_n}$ in $(E /
p^{n+1}E)^\times$ is larger than 1 and divides the prime integer $p'$,
whence it is $p'$.

Next, let $g$ be as above, and let $n \in \N$. By the previous paragraph,
$q_{g+n} = (p')^s q_g$ for some $s$. We now apply the second part of
Proposition \ref{Peuler}, this time setting $x= B^{q_g}, y = 1$. Thus
\[ \nu_p \bigg( \fract{x}{y}{(p')^s} \bigg) = \nu_p((p')^s) = fs \]

\noindent and the result follows.
\end{proof}\hfill

\noindent {\bf Remarks:}
\begin{enumerate}
\item To conclude this discussion, if $d = p_1^{\aaa_1} p_2^{\aaa_2} \dots
p_h^{\aaa_h}$ is coprime to $B$, then to find $D_d$, we find $Q_i =
q_{\aaa_i}(p_i)\ (i=1,2,\dots,h)$. Since $p_1, p_2, \dots, p_h$ all
divide $B^{D_d} - 1$, hence $Q_1, Q_2, \dots, Q_h$ all divide $D_d$.
Moreover, $D_d$ is the least such positive number. Hence
\begin{equation}
D_d = lcm(Q_1, Q_2, \dots, Q_h)
\end{equation}

\item For $E=\Z$, note that $p'=p$ is itself prime, so that $f=1$, and
all our other standing assumptions are satisfied, for any $B>1$.
Moreover, $o_d(B)$ now {\it does} denote the number of digits in the
recurring part of the decimal expansion (in base $B$) of $1/d$, so we can
compute this using the above theorem.\\

\item In the decimal system ($B=10$), note that for the first few primes
$p>5$, we have $g=1$, though $g_3 = 2$ for $p=3$. Thus, $D_3 = D_9 = 1,\
D_{27} = 3,\ D_{81} = 9$, and so on, while $D_7 = 6,\ D_{49} = 6 \times 7
= 42$ etc.

One can ask if this phenomenon actually holds for all primes $p>5$. In
other words, are the orders of $10$ in $(\Z/p\Z)^\times$ and $(\Z /
p^2\Z)^\times$ unequal for all primes $p>5$ ?

The answer is no. The prime $p=487$ satisfies $q_1 = q_2 = p-1=486$ (and
$q_3 = p(p-1)$).

However, is there a base $B$ for which this property does hold (for $E =
\Z$) ?\\
\end{enumerate}

\section{Concluding remarks and questions}\label{Sqns}

We conclude this part with a few questions: all the properties asked
below are known to hold or not to hold for $E=\Z$ or $F[X]$, and we want
to ask whether they hold in general, or if we can characterise all
Euclidean domains with that property.\\

\begin{enumerate}
\item Is there always a ``good" quotient and remainder ? In other words,
can we always find a valuation $\nu$, so that the following holds ?\\

Given $s,0 \neq d \in E$, we can find $q,r$ by the Euclidean algorithm,
such that $s=dq+r,\ r=0$ or $\nu(r) < \nu(d)$, and {\it with the
additional property} that $\nu(dq) \leq \nu(s)$.\\

For example, in $E=\Z$, suppose $s=-11, d=3$. Then $s=(-4)d + 1 = (-3)d +
(-2)$, and we take the good $(q,r)$ to be $(-3,-2)$.\\

\item Does there exist a (sub)multiplicative valuation ? Namely, a
valuation $\nu$ satisfying : if $\nu(a) \leq \nu(b)$ and $\nu(c) \leq
\nu(d)$, then $\nu(ac) \leq \nu(bd)$.\\

\item Does there always exist a valuation satisfying the triangle
inequality ?\\

\item Is there a characterization of all Euclidean domains in which every
prime divides a prime integer (i.e. $p | p'$ in $E$, where $p'$ is prime
in $\Z$, and nonzero in $E$) ?

For example, it is true in $E=\Z$ and in $\Z[i]$ (because $a+bi$ divides
$a^2+b^2 \in \Z$, so every prime divides an integer, and by Lemma
\ref{Lprep9} above, we are done). But this property does not hold in
$E=F[X]$ ($F$ a field), since the polynomial $X-1$ is prime, but does not
divide any integer (or field element).\\
\end{enumerate}

We have a partial answer for this last question (note that we already saw
above that $E$ has to have characteristic zero here).

\begin{prop}
Suppose $(E,\nu)$ is a Euclidean domain but not a field. If the strong
Euclidean algorithm holds, and every prime $p \in E$ divides a prime
integer $p' \in \Z \subset E$, then $E$ is finitely generated as an {\em
$R$-module} over the subring $R = \Z[E^\times]$.
\end{prop}

\begin{proof}
We produce a surjection : $R[X]\to E$, with a {\it monic} polynomial in
the kernel. Note that $E_0 = E^\times \cup \{ 0 \} \subset R$. We first
claim that there is a (set-theoretic) surjection : $E_0[X]
\twoheadrightarrow E$.

Since $E$ is not a field, there exist nonzero nonunits. Pick one with the
least valuation, say $B \in E$, such that $\nu(u) < \nu(B) \leq \nu(y)$
for all nonzero nonunits $y$ and all units $u$ in $E$. By Lemma
\ref{Lprep2}, the evaluation map : $X \mapsto B$ is a surjection from
$E_0[X]$ onto $E$. Hence it extends to an $R$-linear map $: R[X]
\twoheadrightarrow E$.

Next, by factoring $B$ into a product of primes, we can find $n \in \Z$
such that $B|n$. Thus, $Bl = n$ for some $l \in E$. Write $l = p(B)$ for
some polynomial $p$ in $E_0[X]$; we thus get that $n = Bl = Bp(B) =
c_mB^{m+1} + \dots + c_0B$ for some $c_i \in E_0$, and $c_m \neq 0$. Thus
$c_m$ is a unit, and we get that
$$q(B) = c_mB^{m+1} + \dots + c_0 B - n = 0.$$

This means that $c_m^{-1}q(X)$ is in the kernel of the surjection $R[X]
\to E$ (which is the evaluation map at $B$). Since $c_m^{-1}q(X)$ is
monic, hence $1,B,B^2, \dots, B^m$ generate $E$ over $R$, and we are
done.
\end{proof}\hfill

\noindent {\bf Remarks}: The result implies that $E$ is not a polynomial
ring $F[X]$. Moreover, if $E^\times$ is finite, then $E$ is a finitely
generated abelian group. Also, the proof only assumes that the Euclidean
algorithm terminates for all $(s,B)$, where $B$ is any fixed element not
in $E_0$, of least $\nu$-value.\\

Similar questions have been answered; we give a couple of examples. Given
$n \in \N$, we define a Euclidean domain to be {\it of type $n$} if, for
every $s,d \neq 0$ (as in the definition of a Euclidean domain $E$), such
that $d \nmid s$, there exist {\it exactly} $n$ distinct pairs
$(q_i,r_i)$ satisfying the Euclidean algorithm property $s = q_i d + r_i\
\forall 1 \leq i \leq n$. We then have the following two results.

\begin{theorem}[\cite{Ga}]
The only Euclidean domain of type $2$ is $E=\Z$.
\end{theorem}

\begin{theorem}[\cite{Rh,Jo}]
The only Euclidean domains of type $1$ are $E=F$ or $F[X]$, where $F$ is
a field and $X$ is transcendental over $F$. (In particular, $E^\times =
F^\times$.)
\end{theorem}

\noindent Similar results were also shown by Jacobson and Picavet, in
\cite{Ja} and \cite{Pi} respectively.\\

\section*{\bf Part 3 : Quantitative aspects (written with P. Moree)}

According to \cite{Kv}, many people, including Johann Bernoulli III,
C.-F. Gauss, A.H. Beiler, S. Yates, and others, have worked on the
problem of repeating decimals. Repeating decimals were quite a popular
topic of study in the 19th century (cf. Zentralblatt). We also find the
following remarks:\\

All prime numbers coprime to 10 can be divided into three groups:\\

\noindent (1) $p$ such that $D_p = p-1$ (the {\it full-period primes}).\\
(2) $p$ such that $D_p < p-1$ is odd (the {\it odd-period primes}).\\
(3) $p$ such that $D_p < p-1$ is even (the {\it non-full-period
primes}).\\

It was found that the proportion of primes in these three groups had a
relatively stable asymptotic ratio of around 9:8:7, by numerical
computations up to $p=1370471$. We will now explain that it is possible
to be more precise than this.\\

\indent {\it Concerning the primes in} (1). The full-period primes have
been well studied, starting with C.-F. Gauss in his 1801 masterpiece {\it
Disquisitiones Arithmeticae}. Though no written source for this seems to
be available, folklore has it that Gauss conjectured that there are
infinitely many full-period primes.

In September 1927, Emil Artin made a conjecture which implies that the
proportion of primes, modulo which 10 is a primitive root  (i.e. primes
of type (1) in the classification above), should equal a number we now
call {\it Artin's constant} $A = 0.3739558136192 \dots$. More precisely,
we have
$$A = \prod_p \bigg(1 - \frac{1}{p(p-1)}\bigg),$$
where the product is over all primes. Clearly, this is close to 9/24 =
3/8 = 0.375, the ratio mentioned above. Assuming the Generalised Riemann
Hypothesis (GRH) it was proved by Hooley \cite{Ho} that the full-period
primes have proportion $A$.

\indent Artin's original conjecture gave a prediction for the density of
primes $p$ such that a prescribed integer $g$ is a primitive root modulo
$p$. Following numerical calculations by Derrick H. Lehmer and Emma
Lehmer in 1957, Artin corrected his conjecture for certain $g$ (see
\cite{St}), and his corrected conjecture, also attributed to Heilbronn,
was proved by Hooley in the paper cited above, on assuming GRH. Further
generalisations of Artin's primitive root conjecture are discussed, e.g.
in the survey paper \cite{Mo3}.\\

In \cite{He}, Heath-Brown, improving on earlier work by Gupta and Ram
Murty \cite{GR}, proved a result which implies, {\it unconditionally},
that there exist at most two primes $(p > 0)$ - and three squarefree
integers $k > 1$ - for which there are only finitely many full-period
primes (i.e. for which the qualitative version of Artin's primitive root
conjecture fails).\\

{\it Concerning the primes in} (2). Since for an odd prime $p$, the
number $p-1$ is even, the primes $p$ in (2) can be alternatively
described as the primes $p$ such that $D_p$ is odd. These primes also
have been the subject of study (by mathematicians including Sierpinski,
Hasse and Odoni). Without assuming any hypothesis, it can be shown that
the proportion of primes $p$ such that the order of $g$ modulo $p$ with
$g$ any prescribed integer is odd, exists and is a computable rational
number. It turns out  (cf. \cite[Theorem 3.1.3]{Ba} or \cite[Corollary
1]{Mo1}), that if the base $B$ is not of the form $\pm u^2$ or $\pm 2u^2$
for any $u \in \N$, then the proportion of primes $p$ (among all primes)
with odd period, is 1/3. Thus the proportion of the primes in (2) equals
$1/3$.\\

{\it Concerning the primes in} (3). The proportion of primes such that
$D_p$ is even equals $1-1/3=2/3$. Thus, on assuming GRH, the proportion
of primes (3) equals $2/3-A$. As shown in Lemma \ref{Leveniff} above, $B$
has even period mod $p$ with $p\nmid 2B$ iff $p$ divides $B^l+1$ for some
$l \ge 1$. Thus the even period condition modulo $p$ is related to the
divisibility of certain sequences by $p$.\\

\indent In brief, we have arrived at the following result:

\begin{prop}
Assume GRH. Then the natural densities of the sets of primes {\rm (1)},
{\rm (2)}, respectively {\rm (3)}, are given in the table below.
\end{prop}

\medskip
\begin{center}
\begin{tabular}{|c|c|c|c|}\hline
&(1)&(2)&(3)\\
\hline\hline
$\delta$&$A$&${1/3}$&${{2/3}-A}$\\
\hline
$\approx \delta$&$0.37395\ldots$&$0.33333\ldots$&$0.29271\ldots$\\
\hline
{\rm Kvant}&$9/24$&$8/24$&$7/24$\\
\hline
$\approx$&$0.37500\ldots$&$0.33333\ldots$&$0.29166\ldots$\\
\hline
\end{tabular}
\end{center}
\medskip

We next look at the residual index $c_p = [(\Z / p \Z)^\times : \langle B
\rangle]$ (mentioned above), which was also equal to the number of
distinct chains of recurring decimals (in base $B$). Recall that we had
$b_p c_p = p-1$, where $b_p$ is the order of $B$ modulo $p$.

Note that given a fixed $m \in \N$, the number of primes $p$ such that
$b_p = m$ is finite, since there are only finitely many primes less than
$B^m$. We now ask for the density of the set of primes for which the
residual index is fixed, namely $\{ p$ prime : $c_p = m \}$. (Note that
if $m=1$, then this was Artin's conjecture.)

This too has been answered: in \cite{Mu1} we find that assuming the GRH,
the density of this set $N_B^{(m)}$ equals $C_B^{(m)}$ for certain $B$,
which has been explicitly mentioned therein. This density roughly
decreases in the order of $m^{-2}$. Similar results were obtained for
arbitrary $B$ by Wagstaff \cite{Wa}.\\

Finally, we ask how often the order of $B$ modulo $p$ differs from that
modulo $p^2$. Related to this is the question: how often does $p^2$
divide $B^{p-1} - 1$ (or not) ? This is related to the {\it Wieferich
criterion}, which first arose in the study of Fermat's last theorem.

For $x \in \N$, let us denote the primes less than $x$, in each of the
two sets above, by $S_1(x)$ and $S_2(x)$, respectively. Note that every
prime occurs either in $S_1$ or in $S_2$. From \cite{Mu2} we see, roughly
speaking, that the size of $S_2(x)$ can be {\it normally} approximated
(if $B >1$) by a constant times $\ln \ln x$.

Moreover, if we assume the {\it abc}-conjecture to be true, then we infer
(cf. \cite{Si}) that, as $x$ tends to infinity, the cardinality of $\{ p
\leq x : p^2 \nmid (B^{p-1} - 1) \}$ exceeds $c_B \ln x$, where $c_B$
depends on $B$ but not on $x$.\\

\noindent {\bf Acknowledgments :} Thanks are due to Travis Schedler for
finding the counterexample $p=487$ in Section \ref{sec9} above, to
Richard Cudney for first telling me about Artin's conjecture, to
Professor Dan Shapiro, for telling me about the various references
mentioned in Section \ref{Sqns} above, and to the referee for his
suggestions. I also especially thank Professor Pieter Moree for his
remarks and suggestions, that helped bring this paper into its present
form.

\providecommand{\bysame}{\leavevmode\hbox to3em{\hrulefill}\thinspace}
\providecommand{\MR}{\relax\ifhmode\unskip\space\fi MR }
\providecommand{\MRhref}[2]{%
  \href{http://www.ams.org/mathscinet-getitem?mr=#1}{#2}
}
\providecommand{\href}[2]{#2}

\end{document}